\documentclass[12pt]{article}
\usepackage{amsmath, amssymb, amsthm, verbatim,enumerate,enumitem,bbm, mathtools,color,algorithm,float,lineno}
\usepackage[noend]{algpseudocode}
\usepackage{dsfont}
\usepackage[hidelinks]{hyperref}%
\usepackage[capitalize]{cleveref}
\usepackage{tikz-feynman}
\renewcommand{\Pr}[2][]{\mathbb{P}_{#1} \left[ #2 \rule{0mm}{3mm}\right]}

\newtheorem{thm}{Theorem}[section]

\newtheorem{lem}[thm]{Lemma}
\newtheorem{definition}[thm]{Definition}

\newtheorem{claim}[thm]{Claim}
\newtheorem{conj}[thm]{Conjecture}

\makeatletter
\newtheorem*{rep@theorem}{\rep@title}
\newcommand{\newreptheorem}[2]{%
\newenvironment{rep#1}[1]{%
 \def\rep@title{#2 \ref{##1}}%
 \begin{rep@theorem}}%
 {\end{rep@theorem}}}
\makeatother

\newreptheorem{thm}{Theorem}
\newreptheorem{lem}{Lemma}
\newreptheorem{claim}{Claim}

\newcommand{\var}{\varepsilon}

\newcommand{\E}{\mathbb{E}}

\title{Hamilton cycles in regular graphs perturbed by a random 2-factor}
\author{
Cicely (Cece) Henderson
\thanks{Combinatorics and Optimization Department,
University of Waterloo, Waterloo, Ontario N2L 3G1, Canada {\tt c3hender@uwaterloo.ca}.}
\and
Sean Longbrake
\thanks{Department of Mathematics, Emory University, Atlanta, United States {\tt sean.longbrake@emory.edu}}
\and 
Dingjia Mao
\thanks{Department of Mathematics, University of California, Irvine.
Email:  {\tt dingjiam@uci.edu}}
\and 
Patryk Morawski
\thanks{Department of Mathematics, ETH Zürich, Switzerland. This research was conducted while the author was at the Department of Mathematics, University of California, San Diego, United States. Email: \tt{patryk.morawski@ifor.math.ethz.ch}.}}
\date{}

\begin{document}
\maketitle

\begin{abstract}
In this paper, we prove that for each $d \ge 2$, the union of a $d$-regular graph with a uniformly random $2$-factor on the same vertex set is Hamiltonian with high probability. 
This resolves a conjecture by Dragani\'c and Keevash for all values of $d$.
\end{abstract}

\section{Introduction}

A \emph{Hamilton cycle} in a graph $G$ is a cycle that contains all vertices of $G$, and $G$ is \emph{Hamiltonian} if it contains a Hamilton cycle.
Tracing its origins to a problem posed by Hamilton in the 1850s, Hamiltonicity is one of the fundamental notions in graph theory and has been extensively studied (e.g., see the papers \cite{chvatal1972note,draganic2024hamiltonicity,knox2015edge} and surveys \cite{frieze2019hamilton,kuhn2012survey,kuhn2014hamilton}).
In particular, deciding whether a graph contains a Hamilton cycle is among Karp's original list of NP-complete problems \cite{Karp1972}.
It is therefore natural to focus on finding sufficient conditions for a graph to be Hamiltonian. Arguably, the most famous result in this direction is the Dirac's theorem \cite{Dirac1952}, which says that every graph on $n \geq 3$ vertices with minimum degree at least $\delta(G) \geq \frac{n}{2}$ contains a Hamilton cycle.

Dirac's theorem puts a very strong condition on the density of $G$ --- and the constant $\frac{1}{2}$ is tight.
It is natural to ask about the Hamiltonicity of a “typical” graph with a given smaller density.
One example is the binomial random graph $R \sim G(n, p)$ on $n$ vertices, where each edge is included independently and with probability $p$ and where a much lower density suffices to find a Hamilton cycle. 
A result by P\'{o}sa \cite{posa1976hamiltonian}, and independently Korshunov \cite{korshunov1976solution}, states that if $p \gg \frac{n}{\log n}$, then $R \sim G(n, p)$ is Hamiltonian w.h.p.\footnote{With high probability (w.h.p.) means with probability tending to $1$ as $n \to \infty$}.
On the other hand, if $p \ll \frac{\log n}{n}$, then w.h.p. $G$ contains a vertex with degree at most one, so $G$ is not Hamiltonian.
For $d$-regular graphs we need an even lower density.
Already for $d \geq 3$, a uniformly random $d$-regular graph $R$ is Hamiltonian w.h.p. \cite{cooper2002random, krivelevich2001random}, and it is not true for $d =1, 2$.

In \cite{bohman2003many}, Bohman, Frieze and Martin introduced a different random graph model which offers an intermediate perspective between the deterministic and probabilistic settings introduced above.
In the \emph{randomly perturbed graph} model, we start with an arbitrary dense $n$-vertex graph $G_{\alpha}$ with minimum degree at least $\alpha n$ for some constant $\alpha$, and add relatively few random edges to it.
For example, in their paper, Bohman, Frieze and Martin showed that there exists a constant $c = c(\alpha)$ such that the resulting graph $G = G_{\alpha} \cup G(n, c/n)$ is Hamiltonian w.h.p.,  even though $G_{\alpha}$ itself might be far away from being Hamiltonian.
Randomly perturbed graphs $G_{\alpha} \cup G(n,p)$ has been extensively studied since 2003 (e.g. \cite{bottcher2020embedding,joos2020spanning,krivelevich2017bounded}).

What happens if instead of a binomial random graph, we add a uniformly random $d'$-regular graph $R \sim G_{n, d'}$ on top of $G_{\alpha}$?
This question was first studied by Espuny D{\'i}az and Girão \cite{espuny2023hamiltonicity}, who showed that, similarly to the purely random setting, fewer random edges suffice.
Namely, they proved that for $d' = 1$ and any $G_{\alpha}$, the graph $G_{\alpha} \cup G_{n, 1}$ is Hamiltonian w.h.p. provided that $\alpha > \sqrt{2} - 1$, and that the condition on $\alpha$ is tight.
For the case $d' = 2$, they proved that if $\delta(G) > n^{3/4+o(1)}$ then w.h.p. $G \cup G_{n, 2}$ is Hamiltonian.  They also showed that this is not true if $\delta(G) < \frac{1}{5} \log n$.
Draganić and Keevash \cite{draganic2025p} recently improved this bound by showing that if $\delta(G) \geq (1+\var)\sqrt{n \log n /2}$, then $G \cup G_{n, 2}$ is Hamiltonian w.h.p., and that the condition on $\delta(G)$ is tight up to an $\varepsilon$ factor. 

In the same paper, Draganić and Keevash studied the corresponding question in the case when we also require $G$ to be regular.
They show that in this setting already $\delta(G) = \omega(\log^3 n)$ suffices for $G \cup G_{n, 2}$ to be Hamiltonian w.h.p.\footnote{In fact, they show that in this case it is sufficient for $G$ to be only almost regular, i.e., $\Delta(G) = O(\delta(G))$}. They conjecture that this is,  the case for any $d$-regular graph $G$.

\begin{conj}[Conjecture 6.2 in \cite{draganic2025p}]\label{MainConj}
Let $G$ be an $n$-vertex $d$-regular graph, and let $F \sim G_{n,2}$ be a uniformly random $2$-factor on the same vertex set. Then, w.h.p. $G\cup F$ is Hamiltonian.
\end{conj}

In this paper, we resolve this conjecture in the positive for all $d \geq 2$.
\begin{thm}\label{Thm:Main}
    Let $d \geq 2$ be an integer. 
    Let $G$ be an $n$-vertex $d$-regular graph, and let $F \sim G_{n,2}$ be a uniformly random $2$-factor on the same vertex set. Then, w.h.p. $G\cup F$ is Hamiltonian.
\end{thm}

Our proof uses the extension-closure algorithm introduced by Cooper and Frieze \cite{cooper1994hamilton}, which has been used to prove the Hamiltonicity in various sparse random graphs settings (e.g. \cite{cooper19921,cooper1994hamilton,cooper1994hamilton2, frieze1995assignment}). 
We note that, as already observed by Draganić and Keevash, for the case $d = 1$ (i.e., when $G$ is a perfect matching), the union $G\cup F$ is contiguous to the random $3$-regular graph by Theorem 9.34 in \cite{janson2011random}. 
As a consequence of Theorem 9.20 in \cite{janson2011random}, we immediately get that $G\cup F$ is Hamiltonian w.h.p..
Thus \cref{MainConj} is resolved.

The remainder of the paper is organized as follows. In the next subsection, we give an overview of our proof. In \cref{section:prelimiaries}, we collect some preliminary results that we will use in our proof. Our main proof consists of two main phases, where in \cref{section:phase1} we study Phase I via the extension-closure algorithm to eliminate the short cycles, and then in \cref{section:phase2} we study Phase II via the second moment method to close all the cycles into a Hamilton cycle.

\subsection{Proof Overview}\label{section:proof-overview}

To prove \cref{Thm:Main}, we first fix a $d$-regular graph $G$ for $d \geq 2$ and let $F$ be a uniformly random $2$-factor.
We will think of the random process of adding $F$ to $G$ in two steps.
First, we will expose the structure of $F$ --- here, by a standard result (see \cref{lem:permutation properties}), we get that $F$ is a disjoint union of at most $2\log n$ cycles and that most of the cycles are not too short.
In fact, this will be the only randomness about the structure of $F$ we use.
Then, we take a random bijection $\pi: V(G) \to V(F)$ and obtain our random graph $R = \pi(G) \cup F$.

We want to use the second moment method to show that $R$ contains a Hamilton cycle with high probability.
By looking at all possible ways we could connect the cycles of $F$ together using the edges of $\pi(G)$, it is easy to show that the expected number of Hamilton cycles in $R$ is large.
However, the events that a given Hamilton cycle appears in $R$ are heavily dependent on each other --- and the second moment at this stage is unfortunately too large for us.

To circumvent that, we closely follow the lines of Cooper and Frieze \cite{cooper1994hamilton} and split our proof into two parts.
In the first part, we show that using the \emph{extension-closure} algorithm, w.h.p. we can find a $2$-factor in $R$ in which every cycle has length $n_0 = \Omega(\frac{n}{\log n})$. We call this graph $F_{\text{long}}$.
Crucially, during this step we will achieve that while exposing $\pi(v)$ for only a few vertices --- say at most $O(n^{3/4})$.
For the second phase, we will therefore still be able to treat $\pi: V(G) \to V(F)$ essentially as a random permutation.

The extension-closure algorithm works as follows.
We will pick a short cycle $C$ in $F$ and an arbitrary edge $\{u_0, v_0\}$ -- we want to extend this cycle to a cycle of length at least $n_0$.
To that end, we start the following process from $v = v_0$.
First, we expose the vertex $v_G \in V(G)$ with $\pi(v_G) = v$, and expose $w = \pi(w_G)$ for each $w_G \in N_G(v_G)$.
In general, each of these $w$'s will lie on a cycle $C_w$ in $F$ - and will have two neighbors $x_L$ and $x_R$ on this cycle.
In particular, this gives us two paths that start at $v$ and traverse the whole cycle $C_w$ ending at $x_L$ and $x_R$, respectively.

Now, we can repeat this process with the vertices $x_L$ and $x_R$ playing the role of $v$ to find more and longer paths starting at $v$.
We have to be a bit careful to make sure that each $x_L$ is still not exposed and that the vertex $w$ lies outside of the respective path --- we can however show that this happens very rarely.
In particular, we can show that w.h.p. after $O(\log n)$ iterations, we can find a set of at least $\sqrt{n}\log^4 n$ paths starting at $v_0$,  each with a distinct end vertex and each of length at least $n_0$, and that the number of vertices we exposed during this process is small, say $O(n^{3/5})$ (see \cref{lem:tree-growing-lemma}).
We can now repeat the same argument starting from $u_0$ to get a set of at least $\sqrt{n}\log^4n$ long paths starting at $u_0$.
Finally, we show that w.h.p. there will be an edge connecting the endpoints of one of the $v_0$ paths and one of the $u_0$ paths --- which will form our cycle of length at least $n_0$.
Crucially, every other cycle in $F$ will either be completely disjoint to this new cycle or completely lie on it.
We can therefore repeat this process until there are no more short cycles left --- since no new short cycles are ever created, we will need at most $O(\log n)$ iterations. 

In the second step, we aim to use the remaining randomness of $\pi$ to connect the cycles in $F_{\text{long}}$ all together into a Hamilton cycle.
In particular, we consider (possible) Hamilton cycles of the following form:
We pick a set of $O(\log n)$ edges on the cycles in $F_{\text{long}}$ and delete these edges to obtain a collection of $O(\log n)$ paths.
Then, we pick a permutation of these paths - this defines a candidate for a Hamilton cycle obtained by connecting the end vertex of the $i$-th path with the start vertex of the $(i+1)$-th path.
Importantly, each of these Hamilton cycles uses at most $O(\log n)$ random edges --- and the sets of these required random edges for each of the Hamilton cycles considered are mostly disjoint.
With careful calculations, we can therefore show that at this point, the second moment of the random variable counting such Hamilton cycles is small enough to conclude our proof.

\subsection{Notation}

Here we collect notation used throughout our arguments. First, let $G$ be a graph and $x \in V(G)$. We denote the neighborhood of $x$ as $N_G(x)$, or $N(x)$ when the graph is unambiguous. We use $\delta(G)$ and $\Delta(G)$ to refer to the minimum and maximum degree of $G$, respectively. A graph is \emph{$d$-regular} if every vertex has degree $d$. A subgraph is \emph{spanning} if it contains every vertex in the graph. A \emph{$2$-factor} is a 2-regular graph, equivalent to a collection of disjoint cycles spanning the vertex set.  Given a graph $F$, a subset $W \subseteq V(F)$ and a permutation $\pi: V(G) \to V(F)$, we let $F_{\pi, W}$ be the graph on the vertex set $V(F)$ with edges $E(F_{\pi, W}) = E(F) \cup \{\{u, v\} \subseteq W : \{\pi^{-1}(u), \pi^{-1}(v)\} \in E(G) \}$. For any positive integer $n$, we let $[n]:= \{1, 2, ..., n\}$. 

Throughout the proof, we condition on the structure of the uniformly random 2-factor $F$. As we search for a Hamilton cycle, we want certain edges in $G$ to lie between specific vertices in our uniformly random 2-factor, so we use the phrase \emph{expose} to refer to conditioning on the existence of these edges or the lack thereof. Finally, recall that a binomial random variable $\mathrm{Bin}(n, p)$ is the sum of $n$ many independent Bernoulli random variables $I$, where $I = 1$ with probability $p$ and $I = 0$ with probability $1 - p$. 

Our asymptotic notation $O(\cdot), o(\cdot),
\Omega(\cdot), \omega(\cdot)$, etc refers to large $n := |V (G)|$, e.g. $f (n) = \omega(1)$ means $f (n) \to\infty$ as $n \to\infty$. Since all of our calculations are asymptotic, we will often omit floor and ceiling functions whenever they are not crucial. 
Unless explicitly stated otherwise, all logarithms are to the base of $e$.

\section{Preliminaries}\label{section:prelimiaries}

In this section, we collect some preliminary results which we use throughout our proof. First, we frequently use the following estimations to calculate the probability that all random pairs of vertices are edges in a $d$-regular graph $G$.

\begin{lem}\label{obs:probability-of-indicator-variable}
    Let $d \geq 1$, let $G$ be a $d$-regular graph on $n$ vertices and let $W \subseteq V(G)$ be of size $O(n^{3/4})$.
    Let $N$ be a graph on $n - |W|$ vertices and with $m = o(n^{3/4})$ edges, which is a disjoint union of $c$ cycles and some paths.
    Let $\pi: V(G) \setminus W \to V(N)$ be a uniformly random bijection and let $E_N$ be the event that $\{\pi^{-1}(x), \pi^{-1}(y)\} \in E(G)$ for all $\{x, y\} \in E(N)$.
    Then, 
    \[
    \Pr{E_N} \leq (1+o(n^{-1/4}))^{2m} \left(\frac{d}{n}\right)^{m -c}.
    \]
    Moreover, if $N$ is a matching, then 
    \[
    \Pr{E_N} = (1 - o(n^{-1/4}))^{2m} \left(\frac{d}{n} \right)^m.
    \]
    
\end{lem}
\begin{proof}
    Consider the graph $G' = G - W$. Note that $G'$ has at least
    $\frac{dn}{2} - d|W| = \left(1 - o\left(n^{-1/4}\right)\right)\frac{dn}{2}$
    edges. If $N$ is a matching, it follows by inductively conditioning on the success of embedding each edge of $N$ one-by-one.

    Now, let $N$ be a disjoint union of $c$ cycles and some paths. Observe that by removing one edge from each cycle in $N$, we obtain a disjoint union of paths containing $m-c$ total edges. Thus, it suffices to prove the claim for this case.
    Let $v_i^j$ denote the $i$-th vertex on path $j$ in $N$. When conditioning on the images of all previously mapped vertices $\pi(v_{i'}^{j'})$ with $j' < j$, and all $\pi(v_{i''}^j)$ with $i'' < i$, we observe that the probability that $\pi(v_i^j)$ lies in the neighborhood of $\pi(v_{i-1}^j)$ is at most
    $
    \frac{d}{n - |W| - 2m} \leq \left(1 + o\left(n^{-1/4}\right)\right)\frac{d}{n}.
    $
    Proceeding inductively over all non-isolated vertices of $N$ completes the proof.
\end{proof}

The following lemma collects some properties of binomial random variables, where the second one is deduced by the well-known Chernoff's inequality (see, e.g., in \cite{alon2016probabilistic}). 

\begin{lem}\label{lemma:binomial_bounds}
   
    Let $k\geq 1$, and let $B\sim \mathrm{Bin}(n,p)$. 
    Then, $\Pr{B\geq k}\leq \binom{n}{k}p^k$ and $\Pr{B\geq 2np}\leq e^{-np/3}$.
    
\end{lem}

We use a random process generated by some binomial random variables, whose tail bounds can be read as follows.

\begin{lem}\label{lemma:branching_process}
    Let $(S_t)_{t \geq 0}$ be a random process defined by $S_0 =1$ and $S_{t+1} = 3.99S_t - 4B_{t}$, where $(B_t)_{t \geq 0}$ are binomially distributed random variables with $B_t \sim \mathrm{Bin}\left(S_t, \frac{1000}{\log n}\right)$.
    Moreover, let $t_{\max} = \frac{3}{5}\log_4 n$.
    Then, $\Pr{S_{t_{\max}} \geq \sqrt{n}\log^4 n} \geq 1 - O\left(\frac{1}{\log n}\right)$.
\end{lem}

\begin{proof}
    We split the process into two phases.
    In the first phase, when $S_t \leq 1000$, we get by \cref{lemma:binomial_bounds} that with probability at least $1-O(\frac{1}{\log n})$, we have that $S_{t+1} \geq 3.99 S_t$ and thus with probability $1 - O(\frac{1}{\log n})$ we have that $S_t \geq 1000$ for some $t \leq 10$.

    In the second phase, when $S_t\geq 1000$, again by the first part of \cref{lemma:binomial_bounds}, we get that with probability at least $1-O(\frac{1}{\log n})^{S_t/400}\geq 1-O(\frac{1}{(\log n)^2})$, $B_t\leq S_t/400$ and hence $S_{t+1}\geq 3.98S_t$. Therefore, with probability at least $1 - O(\frac{1}{\log n})$, we get that $S_{t_{\max}} \geq \sqrt{n}\log^4 n$.

\end{proof}

The next theorem, due to Vizing \cite{vizing1964estimate}, shows that one can decompose a graph with maximum degree $\Delta$ into at most $\Delta+1$ matchings. In particular, we will use this to reduce a $d$-regular graph $G$  with $d\geq 2$ to a subgraph $G'$ with $2 \le \delta(G') \le \Delta(G') \le 3.$ 

\begin{thm}[Vizing's theorem]\label{thm:Vizing}
Every graph with maximum degree $\Delta$ can be properly edge-colored with $\Delta +1$ colors.
\end{thm}

Finally, we need the following lemma.
It shows that, given two random subsets $X \subseteq V(G)$ and $Y \subseteq V(G)$, if for each $x\in X$ there exists a ``large enough'' set $Y_x\subseteq Y$, then w.h.p. there exists an edge between some $x \in X$ and $y \in Y_x$.

\begin{lem}\label{lem:close-into-cycle}

    Let $G$ be a graph on $n$ vertices with $2\leq \delta(G)\leq \Delta(G)\leq 3$ and let $V'$ be of size $n$.
    Let $W \subseteq V'$ and $W' \subseteq V(G)$ be of sizes $|W| = |W'| = O(n^{3/4})$, and let $\pi: V(G) \setminus W'\to V' \setminus W$ be a uniformly random bijection.
    Let moreover $X \subseteq V' \setminus W$ be of size $\sqrt{n}\log^2 n \leq |X| \leq O(n^{3/4})$, and for each $x \in X$ let $Y_x \subseteq V' \setminus W$ be arbitrary with $|Y_x| \geq \sqrt{n} \log^2 n$.
    Then,  with probability at least $1 - \frac{1}{n^{\Omega(\log n)}}$, there exist $x \in X$ and $y \in Y_x$ such that $\{\pi^{-1}(x), \pi^{-1}(y)\} \in E(G)$. 
\end{lem}

\begin{proof}
    Consider the following process with $\sqrt{n}\log n$ steps.
    In each step, arbitrarily pick an unexposed vertex $x \in X$ such that at least $\frac{1}{2}\sqrt{n} (\log n)^2$ vertices in $Y_x$ are unexposed.
    Then, we expose $x_G = \pi^{-1}(x)$.
    With probability at least $1-O(\frac{1}{n^{1/4}})\geq 1-\frac{\log n}{n^{1/4}}$, the vertex $x_G$ has at least one unexposed neighbor $w_G$ in $G$.
    If this is the case, we expose $w = \pi(w_G)$.
    Then, with probability at least $\frac{\log n}{\sqrt{n}}$ we have $w \in Y_x$, so the probability that a given iteration fails is at most $1 - \frac{\log n}{2 n^{1/2}}$. 

    Now, notice that after $t$ iterations, the number of unexposed vertices in $X$ is at least $|X| - 2t$.
    Similarly, the number of unexposed vertices in $Y_x$ for any $x \in X$ is at least $|Y_x| - 2t$.
    Therefore, we perform at least $\sqrt{n}\log n$ iterations of this process.
    The probability that we fail in every iteration is thus at most
    \[
        \left(1 - \frac{\log n}{4 n^{1/2}}\right)^{\sqrt{n} \log n} \leq e^{-\Omega(\log^2 n)},
    \]
as desired.
\end{proof}

\section{Phase I: Extension-closure algorithm}\label{section:phase1}

In this section, we use the extension-closure algorithm described in \cref{section:proof-overview} to find a 2-factor in $\pi(G)\cup F$ consisting of only long cycles, while exposing few vertices. Our strategy consists of two main steps. In \cref{section:tree-growing}, we start by deleting an edge in $F$ to obtain a ``near-2-factor'', and demonstrate how our paths exploration process works. In \cref{section:closing-cycles}, we apply the paths exploration process in two rounds to close the near-2-factors into cycles and find the desired 2-factor of long cycles. 

\subsection{Finding a large collection of paths starting at a given vertex}\label{section:tree-growing}

Here we describe the process we use to find a large number of paths in our perturbed graph $\pi(G)\cup F$, while only exposing $\pi^{-1}(w)$ for a small number of vertices $w \in V(F)$, as described above.
Before we state our ``paths exploration lemma", we need the following definitions. A \emph{near-$2$-factor} on $n$ vertices is a spanning graph $L$ that is a disjoint union of cycles and exactly one path.
We write $P(L)$ for the unique path in $L$ and write $P(L) \in \mathcal{P}(u, v)$ if the endpoints of this path are $u$ and $v$.  We call a cycle of length at most $n_0\coloneqq  \frac{100 n}{\log n}$ a \emph{short cycle} (and \emph{long cycle} otherwise), and we denote the family of all short cycles in a near-2-factor $L$ by $L_{\text{small}}$.
Using \emph{acceptable rotations}, we grow the path $P(L)$ without creating new small cycles.

\begin{definition}[Acceptable rotation]
    Let $L$ be a near-$2$-factor with $P(L) \in \mathcal{P}(u, v)$, and let $x, w \in V(L)$ with $x \in N_L(w)$.
    A \emph{rotation} $R(L;v,w,x)$ is the near-$2$-factor $L'$ with $E(L') = (E(L) \setminus \{\{w, x\}\}) \cup \{\{v, w\}\}$.
    We say that the rotation $R(L;v,w,x)$ is \emph{acceptable} if 
    \begin{enumerate}[label=(\arabic*)]
        \item $|P(L)| \geq n_0$ implies $|P(L')| \geq n_0$, \emph{and}
        \item any cycle in $L'$ that is not a cycle in $L$ has at least $n_0$ edges.
    \end{enumerate}
\end{definition}
\noindent Notice in particular that we allow the path $P(L')$ to shrink in comparison to $P(L)$ --- but only if both the path $P(L')$ and the new cycle created are longer than $n_0$.

\medskip

We are now ready to state and prove our paths exploration lemma. 
It asserts that, given a near-$2$-factor, we can---with high probability---reach a large collection of near-$2$-factors through acceptable rotations while exposing only a small set of vertices. 
Furthermore, we avoid exposing any vertex from the set $X$ of endpoints of edges reserved on short cycles, and expose only a small fraction of the vertices in the set $X'$ of endpoints of paths starting from the other side.

\begin{lem}[Paths exploration lemma]\label{lem:tree-growing-lemma}
    Let $G$ be a graph on $n$ vertices with $2 \leq \delta(G) \leq \Delta(G) \leq 3$, and let $F$ be a $2$-factor on $n$ vertices.
    Let $W \subseteq V(F)$ with $|W| = O(n^{3/4})$, and let $\pi: V(G) \to V(F)$ be a bijection such that $\pi^{-1}(x)$ is deterministic for all $x \in W$ and the restriction of $\pi$ onto $V(G) \setminus \pi^{-1}(W)$ is a uniformly random bijection into $V(F) \setminus W$.
    Let $X \subseteq V(F) \setminus W$ be of size $|X| \leq 20 \log n$, and let $X' \subseteq V(F) \setminus W$ be of size $|X'|=O(n^{3/4})$.
    Let $v \in V(F) \setminus W$.
    Then, there exists a randomized algorithm that with probability $1 - O(n^{-0.1})$ exposes a set $W' \subseteq V(F)\setminus W$ of size $|W'|=O(n^{3/5})$, and finds a set $S\subseteq V(F)\setminus W$ of size $|S| = O(n^{3/5})$, such that: 
    \begin{enumerate}[label=(T\arabic*)]
        \item\label{tree-growing-property1} $|W' \cap X| = 0$ and $|W' \cap X'| \leq \frac{|X'|}{n^{1/4}}$, \emph{and}
        \item\label{tree-growing-property2} for each near-$2$-factor $L$ in $F_{\pi, W}$ with $P(L) \in \mathcal{P}(u, v)$ for some $u \in V(F)$, with probability at least $1-O(\frac{1}{\log n})$, there is a collection $\mathcal{L}$ of near-$2$-factors of the graph $F_{\pi, W \cup W'}$ with $|\mathcal{L}| \geq \sqrt{n}\log^4 n$, such that for each $L' \in \mathcal{L}$,
        \begin{enumerate}
            \item\label{tree-growing-subproperty1} there exist distinct vertices  $v_{L'}\in S\setminus (W\cup W'\cup X\cup X')$, such that $P(L') \in \mathcal{P}(u, v_{L'})$,
            \item\label{tree-growing-subproperty2} $L'_{\text{small}} \subseteq L_{\text{small}}$, \emph{and}
            \item\label{tree-growing-subproperty3} if $|P(L)| \geq n_0$, then $|P(L')| \geq n_0$.
        \end{enumerate}
    \end{enumerate}
\end{lem}
\begin{proof}
    Throughout the proof, with slight abuse of notation, we allow the set $W'$ to grow implicitly as the algorithm progresses. 
    We define the following subprocedure, \textsc{Traverse}$(\tilde{v}, B)$, which searches for acceptable endpoints starting from an endpoint $\tilde{v}$ of $P(L)$ for some near-$2$-factor $L$, while avoiding a forbidden subset $B \subseteq V(F)$. 
    At the start of this procedure, we assume that $\tilde{v} \notin W \cup W'$, i.e., the vertex $\tilde{v}$ has not yet been exposed. 
    A more intuitive explanation follows the formal definition of the algorithm.

    \begin{algorithm}[H]
\caption{$Traverse(\tilde{v},B)$}\label{procedure:tree-growing}
		
\begin{algorithmic}[1]
\State\label{alg1-expose1}Expose $v_G = \pi^{-1}(\tilde{v})$ and set $W':=W'\cup \{\tilde{v}\}$.
\If{$v_G \in N_G(\pi^{-1}(W \cup W'))$}\label{step:v_in_neighborhood_of_W} 
    \State Return $\emptyset$
\Else
\State\label{alg1-expose2}Let $\{w^{(1)}_G, w^{(2)}_G\} \subseteq N_G(v_G)$ be arbitrary distinct neighbors of $v_G$. 
\State We expose $w^{(j)} = \pi(w_G^{(j)})$ for each $j \in [3]$. Set $W' := W' \cup \{w^{(1)}, w^{(2)}\}$.
    \If{$\{ w^{(1)}, w^{(2)}\} \cap X \neq \emptyset$}\label{step:hit_bad_set}
        \State Terminate with failure.
    \Else  
    \State For each $j \in [2]$, let $N_{F}(w^{(j)}) = \{ x^{(j, 1)}, x^{(j, 2)}\}$.
        \If{Not all $x^{(j, k)}$'s are distinct \emph{OR} for some $j \in [2], k \in [2]$, we have that $x^{(j, k)} \in B \cup W \cup W' \cup X$ }\label{step:checking_xjks}
        \State Return $\emptyset$
        \Else 
        \State Return $\{ x^{(j, k)} : j \in [2], k \in [2] \}$
        \EndIf
    \EndIf
\EndIf 
\end{algorithmic}
\end{algorithm}

The first step of \cref{procedure:tree-growing} is exposing the value $v_G$ of $\pi^{-1}(\tilde{v})$.
Then, we want to expose the value of $w^{(j)} = \pi(w^{(j)}_G)$ for three neighbors $w^{(j)}_G$ of $v_G$ in $G$ --- provided that none of the neighbors of $\pi^{-1}(\tilde{v})$ have been exposed yet.
Now, it might happen that while exposing $w^{(j)}$ we hit a vertex in $X$, in which case we fail to satisfy the conditions of the lemma and terminate with failure.
Otherwise, we look at respectively two neighbors $x^{(j, k)}$ of $w^{(j)}$ on the respective cycles in $F$. Before returning them, we check whether no degeneracies occur.
If at any step we fail, we simply return an empty set.

Before we proceed to the definition of our main algorithm, we analyze a single call of $Traverse(\tilde{v}, B)$.
Notice that since $\tilde{v} \notin W \cup W'$, we get that $\pi^{-1}(\tilde{v})$ is still a uniformly random vertex in $V(G)\setminus(W \cup W')$.
Therefore, the probability that $v_G \in N_G(\pi^{-1}(W \cup W'))$ in line \ref{step:v_in_neighborhood_of_W} is at most $\frac{3|W \cup W'|}{n - |W\cup W'|} = O(n^{-1/4})$, where we used $\Delta(G) \leq 3$ and $|W \cup W'| = O(n^{3/4})$.
Similarly, the probability that $\{ w^{(1)}, w^{(2)}\} \cap X \neq \emptyset$ in line \ref{step:hit_bad_set} is at most $2 \cdot \frac{|X|}{n - |W|} = O(n^{-0.99})$ and the probability that $\{ w^{(1)}, w^{(2)}\} \cap X' \neq \emptyset$ is at most $O(\frac{|X'|}{n})$. 
Finally, since $\Delta(F) = 2$, the probability that in line \ref{step:checking_xjks} we get that $x^{(j, k)} \in B \cup W \cup W' \cup X$ for some $j \in [2]$ and $k \in [2]$ or that not all the vertices $x^{(j, k)}$s are distinct is at most $\frac{4(|B|+|W| + |W'| + |X|)}{n - |W'|} + O(\frac{1}{n}) = O(\frac{|B| + n^{3/4}}{n})$.

   We now set $t_{\max} = \frac{3}{5}\log_4 n$ and claim that by taking $S := S_{t_{\max}}$, the following algorithm satisfies Lemma \ref{lem:tree-growing-lemma}.

    \begin{algorithm}[H]
    
    \caption{Tree growing algorithm $T(L)$}
		\label{procedure:tree-growing-breadth}
        \begin{algorithmic}[1]
        \State Set $S_0 : = \{v\}$.
        \For{$ 1\leq t \leq t_{\max} = \frac{3}{5}\log_4 n$}
        \State Set $S_t : = \emptyset$.
        \For{$x \in S_{t-1}$}
        \If{$x\notin W \cup W'$}
        \State Set $S_{t}\coloneqq S_t\cup Traverse(x,S_{t-1}\cup S_t)$.
        \EndIf
        \EndFor
        \EndFor
        \end{algorithmic}
    \end{algorithm}

    Indeed, each application of algorithm $Traverse(\tilde{v},B)$ outputs at most 4 vertices and exposes at most 3 vertices. So for each $1 \leq t \leq t_{\max}$, we have that $|S_{t}| \leq 4|S_{t-1}|$, and in each iteration $t$ we expose at most $3|S_{t-1}|$ vertices. 
    Thus, we immediately get that $|W'|\leq 3|S_0|+3|S_1|+\ldots+3|S_{t-1}| = O(n^{3/5})$ at the end of one iteration.
    Since we call $Traverse(\tilde{v}, B)$ at most $O(n^{3/5})$ times, we get that $\E[|X \cap W'|] = O(n^{-0.39})$ and $\E[|X' \cap W'|] = O(\frac{|X'|}{n^{2/5}})$.
    In particular, with probability at least $1 - O(n^{-0.1})$ we get that $X \cap W' = \emptyset$ and $|X' \cap W'|\leq \frac{|X'|}{n^{0.1}}$.

    It therefore remains to prove that property \ref{tree-growing-property2} is satisfied.
    To that end, we fix a near-$2$-factor $L$ in $F_{\pi, W}$ with $P(L) \in \mathcal{P}(u, v)$ and define \cref{procedure:acceptable-tree-growing-breadth} to make the analysis easier.
    It will always expose a subset of vertices exposed during the execution of \cref{procedure:tree-growing-breadth} and therefore, we can think of \cref{procedure:acceptable-tree-growing-breadth} as being implicitly executed while running \cref{procedure:tree-growing-breadth}.

    \begin{algorithm}[H]
        \caption{Pruned tree growing algorithm $PT(L)$}
    		\label{procedure:acceptable-tree-growing-breadth}
            \begin{algorithmic}[1]
            \State Set $\mathcal{L}_0 : = \{L\}$.
            \For{$1 \leq t \leq t_{\max} = \frac{3}{5}\log_4 n$}
            \State Set $\mathcal{L}_t \coloneqq \emptyset$. Let $B$ be the set of all non-$u$ endpoints of the paths $P(L')$ for $L' \in \mathcal{L}_{t-1}  $. 
            \For{$L' \in \mathcal{L}_{t-1}$}
                \State Let $v'\in B$ be the non-$u$ endpoint of the path $L'$.
                \If{$v'\notin W\cup W'$}\label{step:already_exposed}\label{step:traverse} 
                        \State Let $\{x^{(j, k)} : j \in [2], k\in[2] \} = Traverse(v', S_{t-1} \cup S_t)$. 
                        \State Let $w^{(1)}$, $w^{(2)}
                        $be the corresponding vertices exposed during $Traverse(v',S_{t-1} \cup S_t)$. \label{step:first_good}
                            \If{all the rotations $L_{j, k}:= R(L'; v', w^{(j)}, x^{(j, k)})$ for $j \in [2], k\in[2]$ are acceptable}
                            \State\label{step:growing_s} Set $\mathcal{L}_t := \mathcal{L}_t \cup \{ L_{j,k} : j \in [2], k \in [2] \}$. 
                            \Else{ Set $\mathcal{L}_t := \mathcal{L}_t$  }
                            \EndIf
            
            \EndIf
            \EndFor
        \EndFor
        \end{algorithmic}
    \end{algorithm}
    
    Before we analyze \cref{procedure:acceptable-tree-growing-breadth}, we first observe that the set $B \subseteq S_t$ represents the set of endpoints that are obtained from acceptable rotations of some specific initial near-2-factor $L$, in contrast to $S_t$, which is the set of all possible endpoints. Also, note that the set of non-$u$ endpoints of $\mathcal{L}_t$ is contained in $S_t$ for all $t$, and thus the set of vertices explored during the execution of Algorithm \ref{procedure:acceptable-tree-growing-breadth} are also explored in the execution of Algorithm \ref{procedure:tree-growing-breadth}.
    Lastly, since at any point, $L'$ is a near-$2$-factor in $F_{\pi, W \cup W'}$ and the vertices $w^{(j)}$ are not exposed, we must have that $N_{L'}(w^{(j)}) = N_F(w^{(j)})$.
    In particular, the rotations $L_{j, k}$ are indeed valid rotations on $L'$. 

    To analyze \cref{procedure:acceptable-tree-growing-breadth}, we first argue that we rarely see $v' \in W \cup W'$ in line \ref{step:already_exposed}.

    \begin{claim}\label{claim:not_many_loops}
        With probability $1 - O(n^{-0.1})$, the following holds.
        In each iteration $t$ we have that $v' \in W\cup W'$ in step \ref{step:already_exposed} of Algorithm \ref{procedure:acceptable-tree-growing-breadth} at most $\frac{|\mathcal{L}_{t-1}|}{n^{0.05}}$ times.
        Moreover, at the end of the process, there are at most $\frac{|\mathcal{L}_{t_{\max}}|}{n^{0.05}}$ near-$2$-factors $L' \in \mathcal{L}_{t_{\max}}$ such that $v'$ is a non-$u$ endpoint of $P(L')$ and $v' \in W \cup W'$.
    \end{claim}
    \begin{proof}[Proof of \cref{claim:not_many_loops}]

        Let $t \geq 1$ and $L \in \mathcal{L}_{t-1}$ with $P(L) \in \mathcal{P}(u, v')$.
        Note that at the time $L$ is added to $\mathcal{L}_{t-1}$ in iteration $t-1$, vertex $v'$ has not been exposed.
        Therefore, if later in iteration $t$ we determine that $v' \in W \cup W'$, it must have been the case that $v' \in \{w^{(1)}, w^{(2)}\}$ during some execution of $Traverse(\tilde{v},B)$.
        In a single execution of $Traverse(\tilde{v},B)$, the probability of this happening is $O(\frac{1}{n})$.
        Moreover, the number of executions of $Traverse(\tilde{v},B)$ between the time we add $L$ into $\mathcal{L}_{t-1}$ in iteration $t-1$ and we check that $v' \in W \cup W'$ in iteration $t$ is at most $|\mathcal{L}_{t-1}|+|\mathcal{L}_t|\leq |S_{t-1}|+|S_t|\leq 5|S_{t-1}|$, since we execute $Traverse(\tilde{v},B)$ exactly once for each $x \in S_{t-1} \cup S_t$.
        This means that the expected number of such vertices $v'$ is at most $O\left(\frac{|\mathcal{L}_{t-1}||S_{t-1}|}{n}\right) = O\left(\frac{|\mathcal{L}_{t-1}|}{n^{2/5}}\right)$.
        Thus, the first claim follows by Markov's inequality and a union bound over all iterations.

        Similarly, for the second claim, notice that at the time we add some $L$ with $P(L) \in \mathcal{P}(u, v')$ to $\mathcal{L}_{t_{\max}}$, we have $v' \notin W \cup W'$.
        After adding $L$ to $\mathcal{L}_{t_{\max}}$, we run $Traverse(\tilde{v},B)$ at most $|S_{t_{\max} -1}| = O(n^{3/5})$ times.
        The second claim thus again follows by Markov's inequality.

    \end{proof}

    Finally, we want to show that $\mathcal{L}_{t_{\max}}$ is large enough.
    The third property follows, since by \cref{claim:not_many_loops} the number of $L' \in \mathcal{L}_{t_{\max}}$ with an exposed end vertex is small w.h.p.
    
    \begin{claim}\label{claim:size-of-L(t')}
        With probability $1 - O(\frac{1}{\log n})$, we have $|\mathcal{L}_{t_{\max}}| \geq 2 \sqrt{n}\log^4 n$.
    \end{claim}

    \begin{proof}[Proof of \cref{claim:size-of-L(t')}]
        By \cref{claim:not_many_loops}, with probability $1 - O(n^{-0.1})$, step \ref{step:traverse} is executed at least $(1 - \frac{1}{n^{0.05}})|\mathcal{L}_{t-1}|$ times during each iteration $t$.
        We condition on this event.

        We now argue that during a given execution of steps \ref{step:first_good}-\ref{step:growing_s} of Algorithm \ref{procedure:acceptable-tree-growing-breadth}, the probability that we add each $L_{j, k}$ into $\mathcal{L}_t$ in step $\ref{step:growing_s}$ is at least $1 - \frac{1000}{\log n}$, independent of the history of the process.
        Indeed, the probability that $Traverse(\tilde{v},B)$ returns an empty set is at most $O(\frac{1}{n^{0.1}})$.
        In the case where one of the rotations in step \ref{step:growing_s} is unacceptable, we selected a vertex $x^{(j,k)}$ among at most $\frac{200n}{\log n}$ vertices that are of distance at most $\frac{100 n}{\log n}$ to the ends of $P(L)$.
        The probability of this occurrence is at most $\frac{200}{\log n}$.
        In particular, we add all six near-2-factors $L_{j, k}$ to $\mathcal{L}_t$ with probability at least $1 - \frac{1000}{\log n}$ as claimed.
 
        Now, in a given iteration $t$, for each execution $g$ of the steps \ref{step:first_good}-\ref{step:growing_s}, let $I_g$ be the indicator random variable for the event that we do not add the sets $L_{j, k}$ to $\mathcal{L}_t$, and let $I = \sum_g I_g$.
        Then, 
        \[
        |\mathcal{L}_t| \geq 4\left(1 - \frac{1}{n^{0.05}}\right)|\mathcal{L}_{t-1}| - 4I \geq 3.99|\mathcal{L}_{t-1}| - 4I.
        \]
        Moreover, $I$ is stochastically dominated\footnote{We say that a random variable $Y$ stochastically dominates a random variable $I$ if $\Pr{Y \geq x} \geq \Pr{I \geq x}$ for all $x \in \mathbb{R}$.} by $\mathrm{Bin}(|\mathcal{L}_{t-1}|, \frac{1000}{\log n})$ and hence $|\mathcal{L}_t|$ stochastically dominates the random variable $3.99|\mathcal{L}_{t-1}| - 4\mathrm{Bin}(|\mathcal{L}_{t-1}|, \frac{1000}{\log n})$.
        The claim now follows by \cref{lemma:branching_process}.
    \end{proof}

    The proof of \cref{lem:tree-growing-lemma} is now completed since (a)--(c) are automatically guaranteed by the algorithms, and the set $W' \subseteq V(F)\setminus W$ of exposed vertices is of size $|W'|=O(n^{3/5})$.
\end{proof}

\subsection{Eliminating short cycles}\label{section:closing-cycles}

In this subsection, we iteratively use \cref{lem:tree-growing-lemma} to eliminate all small cycles in the uniformly random $2$-factor $F$.
As described above, the only randomness about the structure of $F$ we use is the fact that, w.h.p., $F$ does not contain too many short cycles and that $F$ initially contains at least one long cycle.
This follows from a lemma from Dragani{\'c} and Keevash in \cite{draganic2025p} (note that $F\sim G_{n,2}$ can be identified with a uniformly random permutation in $S_n$ conditioned on all cycles having length at least three):

\begin{lem}[Lemma 2.1 in \cite{draganic2025p}]\label{lem:permutation properties}
Let $F \sim G_{n,2}$ be a uniformly random $2$-factor on $n$ vertices. Then for any $\var \in (0,1/2)$, with probability at least $1 - 2\var$ there is a cycle in $F$ of length at least $\var n$. Also, if $t = \omega(1)$ as $n \rightarrow \infty$ then w.h.p. the number of cycles in $F$ of length at most $t$ is $(1+o(1))\log t$; in particular, w.h.p. $F$ contains $(1 + o(1))\log n$ cycles.
\end{lem}

With this in hand, we are ready to state our short cycle elimination lemma.

\begin{lem}[Short cycle elimination lemma]\label{lem:phase1}
    Let $G$ be a $d$-regular graph on $n$ vertices, and let $F$ be a uniformly random $2$-factor on $n$ vertices. 
    Let $\pi:V(G)\to V(F)$ be a uniformly random bijection. 
    Then there exists a randomized algorithm that, with probability at least $1-o(1)$ and by exposing a vertex set $W\subseteq V(G)$ with $|W|\leq O(n^{3/4})$, outputs a $2$-factor in $F_{\pi,W}$ such that each component is a cycle of length at least $n_0=\frac{100n}{\log n}$.
\end{lem}
To prove \cref{lem:phase1}, we eliminate the short cycles one-by-one until none are left.
To remove a particular short cycle $C$, we pick an edge $\{u_0, v_0\}$ on $C$ and delete it from our current $2$-factor.
Then, we use \cref{lem:tree-growing-lemma}, first to find a large collection of paths on the side of $v_0$ and then on the side of $u_0$.
We make one additional intermediate step to ensure the paths are of length at least $n_0$.
Finally, using \cref{lem:close-into-cycle}, we get that w.h.p. there exists an edge closing a long path into a cycle.

\begin{proof}
    We first note that we can instead assume that $2 \leq \delta(G) \leq \Delta(G) \leq 3$. Indeed, if $d\geq 3$, by Vizing's Theorem (\cref{thm:Vizing}), we can partition the edges of $G$ into at most $d+1$ matchings, and taking exactly three of them gives the desired subgraph. 

    First, we consider the structure of $F$.
    By \cref{lem:permutation properties}, we get that w.h.p. $F$ is a union of at most $2 \log n$ cycles with at most $\log \log n$ cycles of length $3$ and at least one cycle of length at least $\frac{n}{\log \log n}$.
    We therefore condition on any such $F$. Recall that a cycle is ``short'' if it has length less than $n_0$. 
    We arbitrarily pick one edge in each cycle of length $3$ and two disjoint edges in each short cycle of length at least $4$.
    Let $X_0$ denote the set of vertices in the edges picked and note that $|X_0| \leq 20\log n$.
    We set $F^{(0)} \coloneqq F$ and $W_0 \coloneqq \emptyset$.

    We now eliminate the short cycles in multiple iterations.
    At the beginning of iteration $t$, we have already exposed a set of vertices $W_{t-1}$ with $W_{t-1}\cap X_{t-1} = \emptyset$ for some $X_{t-1} \subseteq X_0$, and found a $2$-factor $F^{(t-1)}$ in $F_{\pi, W_{t-1}}$ with $F^{(t-1)}_{\text{small}} \subseteq F_{\text{small}}$.
    Assume that $F^{(t-1)}$ contains at least one short cycle $C$.
    We want to expose at most $O(n^{3/5})$ additional vertices to find a $2$-factor $F^{(t)}$ with $F^{(t)}_{\text{small}} \subsetneq F^{(t-1)}_{\text{small}}$ such that $C$ is not a cycle in $F^{(t)}_{\text{small}}$.

    To eliminate $C$ from $F^{(t-1)}$, we pick an edge $\{u_0, v_0\}$ on $C$ with $\{u_0, v_0\} \subseteq X_{t-1}$ and run the algorithm from \cref{lem:tree-growing-lemma} with $v = v_0$, $X = X_{t-1} \setminus \{v_0\}$ and $X'=\emptyset$. 
    With probability at least $1- O(n^{-0.1})$, it exposes a set $B_1$ of vertices with $|B_1| = O(n^{3/5})$ and $B_1 \cap (X_{t-1} \setminus \{v_0\}) = \emptyset$.
    Moreover, with probability $1 - O(\frac{1}{\log n})$, the graph $F_{\pi, W_{t-1} \cup B_1}$ contains a collection $\mathcal{L}_1$ of near-$2$-factors with $|\mathcal{L}_1| = \lfloor\sqrt{n}\log^4 n\rfloor$, such that for each $L \in \mathcal{L}_1$, there exists a distinct $v_L \notin W_{t-1} \cup B_1$ with $P(L) \in \mathcal{P}(u_0, v_L)$ and $L_{\text{small}} \subseteq F^{(t-1)}_{\text{small}}$. 

    We next make a step from $\mathcal{L}_1$ so that all the paths $P(L)$ are of length at least $n_0$.
    To that end, we initially set $\mathcal{L}_2 \coloneqq \emptyset$ and iterate over all $L \in \mathcal{L}_1$.
    If $P(L)$ already has length at least $n_0$, we add $L$ to $\mathcal{L}_2$.
    Otherwise, let $v_L$ be such that $P(L) \in \mathcal{P}(u_0, v_L)$. If $v$ has already been exposed, we discard $L$.
    If not, we expose $v_G \coloneqq\pi^{-1}(v_L)$ and pick an arbitrary $w_G \in N_G(v_G)$ that has not yet been exposed.
    If such a $w_G$ does not exist, we discard $L$, and otherwise we expose $w \coloneqq \pi(w_G)$. 
    With probability $\Omega(\frac{1}{\log n})$, 
    $w$ lies on a cycle of length at least $n_0$ (since there is at least one cycle of length $n_0$ in each $L \in \mathcal{L}_1$) and has an unexposed neighbor $x$ on this cycle, such that $x$ is not an endpoint of any of the paths $P(L)$ for $L \in \mathcal{L}_1 \cup \mathcal{L}_2$. 
    In this event, we add $R(L;v_L,w,x)$ to $\mathcal{L}_2$. Otherwise, we discard $L$.

    Let $B_2$ denote the set of vertices exposed during this step.
    Then we clearly have that $|B_2| \leq 2|\mathcal{L}_1| = O(n^{3/5})$, and thus the probability that $B_2 \cap X_{t-1} \neq \emptyset$ is at most $O(n^{-0.1})$, since $|X_{t - 1}|\leq |X_0| \leq 20 \log n$.
    Moreover, the expected number of times we exposed some $w$ such that $P(L') \in \mathcal{P}(u_0, w)$ for some $L' \in \mathcal{L}_1 \cup \mathcal{L}_2$ is at most $|\mathcal{L}_1| \cdot O(\frac{1}{n^{0.39}})$. In particular, with probability $1 - O(n^{-0.1})$, this does not happen more than $n^{0.45}$ times.
    Therefore, we get that the random variable $|\mathcal{L}_2|$ stochastically dominates $\mathrm{Bin}(\frac{1}{2}\sqrt{n}\log^4n, \Omega(\frac{1}{\log n}))$, and so with probability $1 - O(n^{-0.1})$ we get that $|\mathcal{L}_2| \geq \Omega(\sqrt{n}\log^3 n)$. 
    
    Finally, we let $X'$ be the set of the non-$u_0$ endpoints of the paths $P(L')$ for $L' \in \mathcal{L}_2$ and we run the algorithm from \cref{lem:tree-growing-lemma} again, this time with $v = v_0$, $X = (X_{t-1}\setminus \{u_0, v_0\})$ and $X'$ as just defined.
    With probability $1- O(n^{-0.1})$, this exposes a set $B_3$ of vertices with $|B_3| = O(n^{3/5})$, $B_3 \cap ((X_{t-1}\setminus \{u_0, v_0\}) = \emptyset$ and $|B_3\cap X'| \leq \frac{|X'|}{n^{0.1}}$, and finds a set $S$ of size $O(n^{3/5})$ which is the collection of all non-$v_0$ endpoints of the resulting paths. 
    Moreover, the expected number of $L' \in \mathcal{L}_2$ for which \ref{tree-growing-property2} fails is at most $O(\frac{|\mathcal{L}_2|}{\log n})$, so with probability $1 - O(\frac{1}{\log n})$, there are no more than $\frac{|\mathcal{L}_2|}{2}$ such $L'$'s.
    In particular, with probability $1- O(\frac{1}{\log n})$, for $|\mathcal{L}_2|/3$ unexposed vertices $v_{L'}$ with $L' \in \mathcal{L}_2$, we can find at least $\sqrt{n}\log^3 n$ unexposed and distinct $u_{L'} \in S$, such that $F_{\pi, W_{t-1} \cup B_1 \cup B_2 \cup B_3}$ contains a near-$2$-factor $L''$ with
    \begin{enumerate}
        \item $P(L'') \in \mathcal{P}(v_{L'}, u_{L'})$,
        \item $L''_{\text{small}} \subseteq F^{(t-1)}_{\text{small}}$ and $C$ is not a cycle in $L''$, \emph{and}
        \item $|P(L'')| \geq n_0$.
    \end{enumerate}
    Thus, by \cref{lem:close-into-cycle}, with probability $1 - O(n^{-0.1})$, we can close one of the $P(L'')$'s to a cycle of length at least $n_0$, resulting in a $2$-factor $F^{(t)}$.
    We let $B_4$ be the set of all vertices $v_{L'}$ and the corresponding vertices $u_{L'}$.
    
    At the end of each iteration, we set $W_t = W_{t-1} \cup B_1 \cup B_2 \cup B_3 \cup B_4$ and note that $|B_1 \cup  B_2 \cup B_3 \cup B_4| = O(n^{3/5})$ as desired.
    Moreover, we set $X_t = X_{t-1} \setminus \{ u_0, v_0\}$.

    We repeat this process at most $O(\log n)$ times, so the number of vertices exposed at any point of the algorithm is at most $O(n^{3 / 5} \log n) = O(n^{3/4})$. 
    In particular, we indeed have that all the vertices in $X_{t-1}$ are unexposed at the beginning of iteration $t$.
    The probability that an iteration fails is at most $1 - O(\frac{1}{\log n})$.
    If this happens in an iteration where $C$ has length at least $4$, we rerun the iteration using an alternate edge in $C$ that we chose at the beginning of this process to be in $X_0$. 
    Therefore, the probability that we fail to find a $2$-factor with all cycles being long is at most $\log \log n \cdot O( \frac{1}{\log n}) + 4 \log n \cdot O(\frac{1}{\log^2 n}) + o(1) = o(1)$, as desired.

\end{proof}

\section{Phase II: Connecting the cycles into a Hamilton cycle}\label{section:phase2} 

\begin{proof}[Proof of \cref{Thm:Main}]
Using \cref{lem:phase1}, we get with high probability a 2-factor in $\pi(G) \cup F$ consisting solely of long cycles $C_1, \dots, C_k$, while exposing only $\pi^{-1}(W)$ for a small set $W \subseteq V(F)$ of size $|W| = O(n^{3/4})$. 
As a result, the restriction of $\pi$ to $V(G) \setminus \pi^{-1}(W)$ remains a uniformly random bijection into $V(F) \setminus W$.
In \emph{Phase 2}, we aim to merge the cycles $C_1, \dots, C_k$ using the remaining randomness, ultimately forming a Hamilton cycle in $\pi(G) \cup F$ with high probability. Note that we may assume $k > 1$; otherwise the process is already complete.

To that end, let us now label the vertices in each cycle clockwise as $C_i = \{x_{1, i}, \dots, x_{\ell_i, i}\}$, and let $C'_i = C_i \setminus W$ be the unexposed vertices in each cycle.
We let $c_i = |C_i'|$ and note that $c_i \geq 99 \frac{n}{\log n}$ for all $i \in [k]$ and that $\sum_{i=1}^k c_i \geq 0.9999 \cdot n$.
We set $a = \frac{n}{\log n}$, let $m_i = 2 \lfloor \frac{c_i}{20 \cdot a} \rfloor + 1$ for each $i \in [k]$, and define $m := \sum_{i=1}^{k} m_i$.
We will now define a collection of possible Hamilton cycles that could appear in $F_{\pi, V(F)}$ and later show that one of these Hamilton cycles does indeed appear with high probability.

We consider only Hamilton cycles that can be obtained in the following manner.
On each cycle $C_i$, we pick $m_i$ vertices $x_{j,i} \in C_i'$ at pairwise distance at least $3$ and delete all the edges $\{x_{j - 1, i}, x_{j, i}\}$, where we set $x_{0, i} = x_{\ell_i, i}$.
We let $M$ denote the set of deleted edges and note that there are at least $\prod_{i=1}^k \binom{c_i - 3m_i}{m_i}$ and at most $\prod_{i=1}^k \binom{c_i}{m_i}$ possible choices for $M$.

Next, we temporarily relabel the edges in $M$ as follows.
On each cycle $C_i$, we identify the smallest $j_i$ such that $\{x_{j_i - 1, i} , x_{j_i, i}\} \in M$.
We set $u_1 = x_{j_1 - 1, 1}$ and $ v_1 = x_{j_1, 1}$, and then continue clockwise around $C_1$ defining $u_2, v_2, \dots, u_{m_1}, v_{m_1} \in C_1'$ in order.
We then set $u_{m_1 + 1} = x_{j_2 - 1, 2}$ and $v_{m_1 + 1} = x_{j_2, 2}$, and so on. 
After removing $M$ from the cycles $C_1, \dots, C_k$, we are left with $m$ path sections $P_j \in \mathcal{P}(u_{\phi_M(j)}, v_j)$ in $F_{\pi, W}$ for some permutation $\phi_M$.
We note that $\phi_M$ is uniquely determined by the choice of $M$ and that since $m_i$ is odd, we get that $\phi_M$ is an even permutation for each $i \in [k]$.

We want to rejoin these paths into a Hamilton cycle using the (random) edges of $G$.
To that end, let $H_m$ denote the set of all cyclic permutations on $[m]$.
Then, each $\rho \in H_m$ defines a distinct Hamilton cycle obtained by joining $P_i$ to $P_{\rho(i)}$ using the edge $\{v_i, u_{\phi_M(\rho(i))}\}$.
To obtain a reduction in the second moment we restrict the choices of $\rho$ to $\rho \in R_M = \{\rho \in H_m : \phi_M \circ \rho \in H_m \}$.

For each possible choice of $M$ and $\rho \in R_M$ we now let $H_{M, \rho}$ be the indicator random variable for the event that $\{\pi^{-1}(v_i), \pi^{-1}(u_{\phi_M(\rho(i))} )\} \in E(G)$ for all $i \in [m]$, where $u_i$ and $v_i$ are implicitly defined by $M$ as above.
Moreover, let $H = \sum_{M, \rho} H_{M, \rho}$, where the sum is over all possible choices of $M$ and $\rho \in R_M$.
We claim that
\begin{claim}\label{claim:first_moment}
    $\mathbb{E}[H] = \omega(1)$ and $\mathbb{E}[H^2] = (1+o(1))\mathbb{E}[H]^2$.
\end{claim}
We defer the proof to the subsequent subsection.
Assuming \cref{claim:first_moment}  holds, we immediately get \cref{Thm:Main}, since
\[
    \Pr{F_{\pi, V(F)} \text{ is not Hamiltonian}} \leq \Pr{H = 0} \leq 1 - \frac{\mathbb{E}[H]^2}{\mathbb{E}[H^2]} = o(1).
\]
\end{proof}
\subsection{Moments computations}
\begin{proof}[Proof of \cref{claim:first_moment}]
In this section, we want to bound $\mathbb{E}[H]$ and $\mathbb{E}[H^2]$ to prove \cref{claim:first_moment}.
To that end, we first notice that
\[
    0.099 \log n \leq \frac{n}{10 \cdot a} - 4k \leq m \leq \frac{n}{10 \cdot a} + 4k \leq 0.101 \log n.
\]
Moreover, 
\[
    k \leq \frac{\log n}{100}, \qquad m_i \geq 5 \qquad \text{and} \qquad \frac{c_i}{m_i} \geq 9  a.
\]
By Lemma 3 in \cite{frieze2001hamilton} we also get that $(m-2)! \leq |R_M| \leq (m-1)!$ for any $M$ chosen as above.
Now, by \cref{obs:probability-of-indicator-variable} we get that for each $M$ and $\rho \in R_M$,
\[
    \mathbb{E}[H_{M, \rho}] = (1+o(1)) \left( \frac{d}{n} \right)^m,
\]
and therefore
\begin{align*}
    \mathbb{E}[H] &= \sum_{H, \rho} \mathbb{E}[H_{M, \rho}] \\
    &\geq \frac{1}{2} \sum_{M, \rho} \left(\frac{d}{n}\right)^m\\
    &\geq \frac{1}{2} \left(\frac{d}{n} \right)^m (m-2)! \prod_{i=1}^k \binom{c_i - 3m_i}{m_i}\\
    &\geq \frac{1}{2}\left(\frac{d}{n} \right)^m \left(\frac{m-2}{e}\right)^{m-2} \prod_{i=1}^k \frac{1}{2\pi m_i} \left( \frac{c_ie}{m_i} \right)^{m_i}\\
    &\geq \frac{1}{2m^2 \cdot (2\pi)^k} \left( \frac{dm}{en} \cdot 8 \cdot ea \right)^m \\
    &\geq n^{-0.02} \left( \frac{d}{1.5} \right)^m,
\end{align*}
which proves $\mathbb{E}[H] = \omega(1)$ since $d \geq 2$.

We now want to bound $\mathbb{E}[H^2]$ --- to do that we need to understand $\mathbb{E}[H_{M, \rho} \cdot H_{M', \rho'}]$ for distinct pairs $(M, \rho)$ and $(M', \rho')$.
We therefore fix one such pair and let $N$ be the set of edges connecting the path segments in $(M, \rho)$ , i.e., $N = \{ \{v_i, u_{\phi_M(\rho(i))} \}: i \in [m]\}$, where $v_i, u_i$ are defined implicitly by $(M, \rho)$.
We define $N'$ analogously for $(M', \rho')$.

We first note that both $N$ and $N'$ are matchings.
Moreover, if $N \cap N' = \emptyset$ and $N \cup N'$ is acyclic, then by \cref{obs:probability-of-indicator-variable}, we get that $\mathbb{E}[H_{M, \rho} \cdot H_{M', \rho'}] \leq (1+o(1))\mathbb{E}[H_{M, \rho}] \mathbb{E}[H_{M', \rho'}]$.
In general, however, it is possible to have overlapping edges between $N$ and $N'$, as well as cycles in $N \cup N'$.

Therefore, let $t = |N \cap N'|$, let $c$ be the number of cycles in $N \cup N'$, and let $\ell \geq 2c$ be the number of $N$ edges in all cycles in $N \cup N'$.
By \cref{obs:probability-of-indicator-variable} we get that
\[
    \mathbb{E}[H_{M, \rho} \cdot H_{M', \rho'}] \leq \left( \frac{d}{n} \right)^{2m -t - c} \leq (1 + o(1)) \left( \frac{n}{d} \right)^{t +c} \mathbb{E}[H_{M, \rho}] \mathbb{E}[  H_{M', \rho'}].
\]

For a fixed $(M, \rho)$ we now aim to bound the number of $(M', \rho')$ giving rise to specific values of $t$, $c$, and $\ell$.
This requires the following observation.
Let $s = |\{(e, f) \in M \times M': e\cap f \neq \emptyset\}|$ and note first that by construction we get $t + \ell \leq s$.
Moreover, using that $\rho \in R_M$, we now argue that if $(M, \rho)$ and $(M', \rho')$ are distinct then we always have $t < s$.

Indeed, suppose that $s = t > 0$ and let $\lambda = \phi_M \circ \rho$.
We define an auxiliary graph $C$ on the vertex set $S = \{(e, f) \in M \times M' : e \cap f \neq \emptyset \}$, where $\{(a,b), (c, d)\} \in E(C)$ if and only if there exists an edge $g \in N \cap N'$ with one endpoint in $a \cap b$ and the other endpoint in $c \cap d$. Since edges in $M'$ are at pairwise distance at least $2$, for each $e \in M$ there is at most one $f \in M'$ with $e \cap f \neq \emptyset$.
Moreover, since $N \cap N'$ is a matching, the auxiliary graph $C$ has maximum degree at most $2$.
Additionally, a vertex $(e, f) \in S$ can have degree $2$ in $C$ only if $e = f$.

Now, by assumption $|E(C)| = t = s = |S|$, so the graph $C$ is $2$-regular. We now pick any $\{v_i, u_{\lambda(i)}\} \in N \cap N'$ and notice that we immediately get $\{u_i, v_i\}, \{u_{\lambda(i)}, v_{\lambda(i)} \} \in M \cap M'$.
Moreover, since $C$ is $2$-regular, we get that $\{ v_{\lambda(i)}, u_{\lambda^2(i)}\} \in N \cap N'$.
Repeating this argument inductively, we get that $\{u_{\lambda^k(i)}, v_{\lambda^{k}(i)}\} \in M \cap M'$ for all $k \geq 0$.
Since $\rho \in R_M$, i.e., $\lambda$ is cyclic, we get that $M = M'$ as claimed.



We now count the number of $(M', \rho')$ that give rise to specific values of $(s, t, c, \ell)$ for fixed $(M,\rho)$ in the following way.
We first pick $s$ edges of $M$ and for each of them pick one of three possible choices for the edge in $M'$ that contributes towards $s$.
Then, we pick $\ell$ of these $s$ edges to witness the $N$-edges in the cycles in $N \cup N'$ and $t$ of the remaining $s - \ell$ edges to witness the edges $N \cap N'$.
Finally, we choose the way these edges give rise to $c$ cycles, pick the remaining edges for $M'$ and fix the remainder of $\rho'$.
By this logic, for fixed $(M, \rho)$, the number of $(M', \rho')$ giving rise to specific $(s, t, c, \ell)$ is at most
\begin{align*}
    &\binom{m}{s} \cdot 3^s \cdot \binom{s - \ell}{t}\binom{s}{\ell} \cdot  \ell^c \cdot \max_{\sigma_1 + \dots + \sigma_k = s} \left[ \prod_{i = 1}^k \binom{c_i}{m_i - \sigma_i}  \right] (m-t -\ell)! \\
    \leq & 2^km! \cdot \left[ \prod_{i = 1}^k \binom{c_i - 3m_i}{m_i}  \right] \cdot \frac{m^{s - \ell - t}}{s!} \binom{s}{t} s^{\ell} \cdot  \ell^c  \cdot (3a)^{-s},
\end{align*}
where we used that
\[
    \frac{\binom{c_i}{m_i - \sigma_i}}{\binom{c_i - 3m_i}{m_i}} \leq 2 \left( \frac{m_i}{c_i}\right)^{\sigma_i} \leq 2 \left( 9a \right)^{- \sigma_i}.
\]
Combining the above observations, still for a fixed $(M, \rho)$, we therefore get that
\begin{align*}
    \frac{\sum_{M', \rho'} \mathbb{E}[H_{M, \rho} \cdot H_{M'. \rho'}]}{\sum_{M', \rho'} \mathbb{E}[H_{M, \rho}] \mathbb{E}[ H_{M'. \rho'}]} &\leq 1 + o(1) + 2^k m^2\sum_{s, t, c, \ell} \frac{m^{s - \ell - t}}{s!} \binom{s}{t} s^{\ell} \cdot  \ell^c  \cdot (3a)^{-s} \cdot \left( \frac{n}{d} \right)^{t+c}.
\end{align*}
For the cases where $\ell = c = 0$, using $t \leq s-1$, we can upper bound each summand by  
\begin{align*}
    \frac{1}{s!} \left( \frac{m}{3a} \right)^{s} \binom{s}{t} \left( \frac{n}{dm} \right)^{t} \leq \frac{1}{s!} \left( \frac{m}{3a} \right)^{s} \cdot s \cdot \left( \frac{n}{dm} + 1 \right)^{s-1} \leq \frac{dm^3}{n}\cdot \frac{1}{s!} \cdot \left( \frac{n}{2da} \right)^s,
\end{align*}
where we used that $\binom{s }{t}x^t \leq s\sum_{i = 1}^{s - 1}\binom{s -1}{i}x^i = s(1 + x)^{s - 1}$. For the remaining cases, where $\ell \geq 2c \geq 2$, and $t \leq s - \ell$, we similarly upper bound each summand by
\begin{align*}
    \frac{m^{s-\ell}}{s!} \binom{s-\ell}{t} \left( \frac{n}{dm} \right)^t \cdot s^{\ell} \cdot \ell^c \left( \frac{n}{d} \right)^c (3a)^{-s} &\leq \frac{m^{s-\ell}}{s!} \left( \frac{n}{dm} + 1 \right)^{s-\ell} \left(s^2 \cdot \ell^2 \right)^{\ell/2} (3a)^{-s}\cdot \left(\frac{n}{d}\right)^{\ell / 2}\\
    &\leq \frac{1}{s!} \left( \frac{n}{2ad} \right)^s \left( \frac{ds^2\ell^2}{n} \right)^{\ell /2}\\
    &\leq \frac{dm^4}{n}\cdot \frac{1}{s!} \cdot \left( \frac{n}{2ad} \right)^s.
\end{align*}
In particular, as $s \leq m$ and so the number of choices for $(t, c, \ell)$ is at most $m^3$, we get that
\begin{align*}
     \frac{\sum_{M', \rho'} \mathbb{E}[H_{M, \rho} \cdot H_{M'. \rho'}]}{\sum_{M', \rho'} \mathbb{E}[H_{M, \rho}] \mathbb{E}[ H_{M'. \rho'}]} &\leq 1 + o(1) + n^{-0.9} \cdot \sum_{s=1}^{m} \frac{1}{s!} \left( \frac{n}{2da} \right)^s \\
     &\leq 1+ o(1) +n^{-0.9} \cdot \sum_{s=1}^m \left( \frac{e \log n}{4s}\right)^s \\
     &\leq 1+ o(1),
\end{align*}
since $(\frac{e \log n}{4s})^s$ is maximized when $s = \log n / 4$ and can therefore be upper bounded by $n^{0.3}$.
Consequently,
\[
    \frac{\mathbb{E}[H^2]}{\mathbb{E}[H]^2} =\frac{\sum_{M, \rho} \sum_{M', \rho'} \mathbb{E}[H_{M, \rho} \cdot H_{M'. \rho'}]}{\sum_{M, \rho}\sum_{M', \rho'} \mathbb{E}[H_{M, \rho}] \mathbb{E}[ H_{M'. \rho'}]} = 1 + o(1),
\]
which proves $\mathbb{E}[H^2] = (1+o(1))\mathbb{E}[H]^2$. This completes the proof of \Cref{claim:first_moment}.

\end{proof}

\section{Concluding remarks}

As mentioned in the introduction, Draganić and Keevash~\cite{draganic2025p} show that for $d=\omega(\log^{3} n)$, it suffices that $G$ is almost $d$-regular, i.e., $\delta(G)\ge d$ and $\Delta(G)=O(\delta(G))$. We expect that our proof can be adapted to work for almost $d$-regular $G$ for all $d\ge 2$. In fact, we believe one can weaken the assumption further to $\delta(G)\ge 2$ and $\Delta(G)=o(\delta(G)\log n)$ and still obtain that $G\cup H$ is Hamiltonian w.h.p. For simplicity of presentation, we do not pursue this stronger statement here.
We would nonetheless find it interesting to determine whether the assumption on $G$ can be relaxed even further.



\section*{Acknowledgements}
We would like to thank the organizers and the participants of the Desert Discrete Math Workshop II (DDMWII) supported by NSF CAREER Grant DMS-2146406, where this research project was initiated and a part of it was completed. 
Special thanks to Asaf Ferber and Marcelo Sales for organizing the workshop and introducing us to the problem. 
We would also like to thank Lior Gishboliner for constructive and detailed comments on an earlier draft of this paper. Additionally, we thank Asaf Ferber, Luke Postle, and Liana Yepremyan for proofreading our draft and giving us constructive criticism.

\bibliographystyle{abbrv}
	\bibliography{references}
\end{document}